\documentclass[11 pt, twoside, reqno]{amsart}
\usepackage{parskip}
\usepackage{amsthm,amsmath,amssymb,oldgerm}
\usepackage[a4paper, margin=2.5cm]{geometry}
\usepackage{mathrsfs}
\usepackage{verbatim}
\usepackage{hyperref}
\usepackage{pdfpages}
\theoremstyle{plain}
\newtheorem{thm}{Theorem}[section]
\newtheorem*{theorem*}{Main theorem}
\newtheorem{mainthm}[]{Main Theorem}[]
\newtheorem{maincor}[mainthm]{Corollary}

\newtheorem{prop}[thm]{Proposition}
\newtheorem{cor}[thm]{Corollary}
\theoremstyle{definition}
\theoremstyle{definition}
\newtheorem{rem}[thm]{Remark}

\title[Estimates and asymptotics of Teichm\"uller modular forms]{Estimates and asymptotics of Teichm\"uller modular forms}
\author{Anilatmaja Aryasomayajula}
\address{Department of Mathematics, Indian Institute of Science Education and Research (IISER) Tirupati, 
Transit campus at Sri Rama Engineering College, Karkambadi Road,
Mangalam (B.O),Tirupati-517507, India.}
\email{anil.arya@iisertirupati.ac.in}
\author{Debasish Sadhukhan}
\address{Department of Mathematics, Indian Institute of Science Education and Research (IISER) Tirupati, 
Transit campus at Sri Rama Engineering College, Karkambadi Road,
Mangalam (B.O),Tirupati-517507, India.}
\email{debasish@students.iisertirupati.ac.in}
\date{\today}
\subjclass[2020]{}
\begin{document}
\begin{abstract}
In this article, we derive estimates of Teichm\"uller modular forms, and associated invariants. Let $\mathcal{M}_{g}$ denote the moduli space of compact hyperbolic Riemann surfaces of genus $g\geq 2$, and let $\overline{M}_{g}$ be the Deligne-Mumford compactification  of $\mathcal{M}_{g}$, and we denote its boundary by $\partial\mathcal{M}_{g}$. Let $\pi:\mathcal{C}_{g}\longrightarrow\mathcal{M}_{g}$ be the universal surface.  For any $n\geq 1$, let $\Lambda_{n}:=\pi_{\ast}(T_{v}\mathcal{C}_{g})^{n}$, where $T_{v}\mathcal{C}_{g}$ denotes the vertical holomorphic tangent bundle of the fibration $\pi$, and the fiber of $\Lambda_{n}$ over any $X\in\mathcal{M}_{g}$ is equal to $H^{0}(X,\Omega_{X}^{\otimes n})$, the space of holomorphic differentials of degree-$n$, defined over the Riemann surface $X$. Let $\lambda_{n}:=\mathrm{det}(\Lambda_{n})$ denote the determinant line bundle of the vector bundle $\Lambda_{n}$, whose sections are known as Teichm\"uller modular forms. The complex vector space of Teichm\"uller modular forms is equipped with Quillen metric, which is denoted by $\|\cdot\|_{\mathrm{Qu}}$. 

\vspace{0.1cm}\noindent 
Let $\lbrace X_t\rbrace$ denote a family of compact hyperbolic Riemann surfaces on $\mathcal{M}_{g}$, and as $t$ approaches zero, $X_t$ degenerates to a noncompact hyperbolic Riemann surface $M_0\in \partial\mathcal{M}_g $, which is of finite hyperbolic volume. In this article, we derive estimates of the Quillen metric of a certain family of Teichm\"uller modular forms defined over the family of compact hyperbolic Riemann surfaces $\lbrace X_t\rbrace$, and ascertain their behaviour, as $t$ approaches zero. Using these estimates, we also derive certain crude estimates of Mumford forms. 
\end{abstract}
\maketitle
\section{Introduction}
\subsection{History}\label{sec-1.1}
Estimates of sections of holomorphic line bundles defined over complex manifolds is an active area of research in complex geometry. Especially, estimates of automorphic forms, which can be viewed as sections of a holomorphic line bundle defined over a Shimura variety, are of great interest, both in complex geometry, and in number theory. Furthermore, estimates in the setting of noncompact complex manifolds are difficult to obtain. 

\vspace{0.1cm}
Through Polyakov's quantization of string theory in \cite{poly1} and \cite{poly2}, followed by Belavin and Kniznik's formula for the partition function in the setting of Bosonic string theory (\cite{bk}), the underlying connection between sections of certain holomorphic line bundles defined over moduli space of Rieamann surfaces, and the partition function in Bosonic string  theory are well documented. These sections are known as Teichm\"uller modular forms, which are the natural analogues of classical modular forms and Siegel modular forms, which are defined over  modular curves and Siegel modular varieties, respectively. 

\vspace{0.1cm}
As sections of a holomorphic line bundle, classical modular forms and Siegel modular forms are linked via the Torelli map. Estimates of classical cusp forms and Siegel cusp forms, especially, sub-convexity estimates, are of great interest in number theory. 

\vspace{0.1cm}
In this article, we first obtain lower and upper bounds for the Selberg-zeta function, using which we derive estimates of Teichm\"uller modular forms, and eventually derive estimates of  Mumford forms. Our estimates are not optimal in certain aspects, but they can be improved by combining with results from \cite{MT}, which will be dealt in a future article. 

\vspace{0.1cm}
\subsection{Statement of results}\label{sec-1.2}
Let $\mathcal{M}_{g}$ denote the moduli space of compact hyperbolic Riemann surfaces of genus $g\geq 2$, which is a complex orbifold of dimension $3g-3$. Let $\overline{\mathcal{M}}_{g}$ denote the Deligne-Mumford compactification  of $\mathcal{M}_{g}$, and let $\partial \mathcal{M}_{g}:=\overline{\mathcal{M}}_{g}\backslash \mathcal{M}_{g}$ denote the boundary of $\overline{\mathcal{M}}_{g}$. 

\vspace{0.1cm}
Let $\pi:\mathcal{C}_{g}\longrightarrow \mathcal{M}_{g}$ denote the universal surface. Let $T_{v}\mathcal{C}_{g}$ denotes the vertical holomorphic tangent bundle of the fibration 
$\pi$. For any $n\geq $1, put $\Lambda_{n}:=\pi_{\ast}(T_{v}\mathcal{C}_{g})^{n}$, and the fiber of $\Lambda_{n}$ over any $X\in\mathcal{M}_{g}$ is equal to $H^{0}(X,\Omega_{X}^{\otimes n})$, the space of holomorphic differentials of degree-$n$, defined over the Riemann surface $X$. Let $H^{0}(\mathcal{M}_{g},\lambda_{n})$ denote the space of global holomorphic sections of $\lambda_{n}:=\mathrm{det}(\Lambda_{n})$, the determinant line bundle of the vector bundle $\Lambda_{n}$. Global sections of the line bundle $\lambda_n$ are known as Teichm\"uller modular forms. 

\vspace{0.1cm}
The complex vector space $H^{0}(\mathcal{M}_{g},\lambda_{n})$ is equipped with the Quillen metric $\|\cdot\|_{\mathrm{Qu}}$. Let $\varphi^n\in \lambda_{n}$ be any section with trivialization $\varphi^n:=\varphi_{1}\wedge\cdots\wedge \varphi_{g_{n}}$ at  $X\in \mathcal{M}_{g}$, where $\lbrace \varphi_{1},\ldots,\varphi_{g_{n}}\rbrace$ denotes a fixed basis of $H^{0}(X, \Omega_{X}^{\otimes n})$ (which is as described in section \ref{sec-2.4}), and $g_{n}:=(2n-1)(g-1)+\delta_{1}(n)$ denotes the dimension of the complex vector space $H^{0}(X, \Omega_{X}^{\otimes n})$, which is equipped with the $L^2$ inner-product $\langle\cdot,\cdot\rangle_{\mathrm{hyp}}$. Then, the Quillen metric of $\varphi^n$ at a $X\in\mathcal{M}_{g}$, is given by the following formula
\begin{align*}
\|\varphi^n(X)\|_{\mathrm{Qu}}^2:=\frac{\big|\mathrm{det}N_{\varphi^n}(X)\big|}{\big|\mathrm{det}^{\ast}\Delta_{\mathrm{hyp},n} (X)\big|},
\end{align*}
where $N_{\varphi^n}(X)$ denotes the matrix
\begin{align*}
N_{\varphi^n}(X):=\big(\left\langle\varphi_{i},\varphi_{j}\right\rangle_{\mathrm{hyp},n}\big)_{1\leq i,j\leq g_{n}},
\end{align*}
and the inner-product $\langle\cdot,\cdot\rangle_{\mathrm{hyp}}$ is as defined in equation \eqref{hyp-met}, and $\mathrm{det}^{\ast}\Delta_{\mathrm{hyp},n} (X)$ denotes the regularized determinant of the hyperbolic Laplacian acting on smooth $n$-differentials, defined over $X$. 

\vspace{0.1cm}
For $g,n\geq 1$, from \cite{sar}, \cite{DP}, and \cite{giovanni}, we have 
\begin{align*}
\mathrm{det}^{\ast}\Delta_{\mathrm{hyp},n} (X)=C_{g,n}Z_{X}(n),
\end{align*}
where $Z_{X}(n)$ denotes the Selberg-zeta function, and $C_{g,n}$ is a constant which only depends on $g$ and $n$. The constant has been explicitly computed in \cite{giovanni}, and it is as described in equation \eqref{cgn2}.  

\vspace{0.1cm}
Furthermore, from Mumford isomorphism 
\begin{align*}
H^{0}\big(\mathcal{M}_{g},\lambda_{n}\big)\simeq H^{0}\big(\mathcal{M}_{g}, \lambda_1^{\otimes (6n^2-6n+1)}\big),
\end{align*}
which is an isometry of vector spaces, we infer that the line bundle $\lambda_{n}\otimes \lambda_{1}^{-\otimes(6n^{2}-6n+1)}$ is trivial. Mumford showed that the line bundle $\lambda_{n}\otimes \lambda_{1}^{-\otimes(6n^{2}-6n+1)}$ admits a non vanishing, non constant holomorphic section. This section is unique up to a constant, and is called the Mumford form, which we denote by $\mu_{g,n}$. For $n=2$, Mumford forms appear in the Polyakov formula for the partition function in Bosonic string theory.  

\vspace{0.1cm}
Let $\lbrace X_{t}\rbrace$ be a family of compact hyperbolic Riemann surfaces of genus $g\geq 2$ on $M_{g}$, which is as described in section \ref{sec-2.4}. As $t$ approaches zero, $X_t$  degenerates to a noncompact hyperbolic Riemann surface $X_{0}$ of finite hyperbolic volume with $j_{0}$ punctures. 

\vspace{0.1cm}
We use the following notation. For any two real-valued functions $f$ and $g$,  $f\ll_{M_{0}} g$ implies there exists a constant $c$, which depends only on $M_{0}$ such that $f\leq cg$.

\vspace{0.1cm}
The first main result of the article is the following theorem, where we estimate the Selberg-zeta function, and hence, $\mathrm{det}^{\ast}\Delta_{\mathrm{hyp},n}$. 

\vspace{0.1cm}
\begin{mainthm}\label{mainthm1}
With notation as above, for any $n\geq 2$ and $X\in \lbrace X_{t}\rbrace$, we have the following inequalities
\begin{align}\label{mainthm1-est}
&Z_{X}(2)\ll Z_{X}(n)\ll(4n^2-4n-3)Z_{X}(2)\bigg(\prod_{j=1}^{j_{0}}e^{c(g,\ell_{j})/\ell_{j}^{2}}\bigg),\notag\\&
\mathrm{where,\,\,for}\,\,1\leq j\leq j_{0},\,\,c(g,\ell_{j}):=e^{\frac{160\pi(g-1)}{\ell_j}},
\end{align}
where $\ell_{j}$ denotes the length of the closed geodesics, which approach zero, as $X$ approaches $X_{0}$.
\end{mainthm}

\vspace{0.1cm}
Observe that Main Theorem \ref{mainthm1} is a result at the level of Riemann surfaces, and number-theoretic in flavor. However, it will be useful in proving Corollaries \ref{mainthm2} and \ref{mainthm3}, which are geometric in nature, with plausible repercussions in string theory. Furthermore, in \cite{wo}, Wolpert proved that
\begin{align}\label{woconj}
Z_{X}(n)=O_{K}\big(Z_{X}(2)\big),
\end{align} 
where the implied constant depends on 
\begin{align*}
n\in K\subset \big\lbrace s\in\mathbb{C}\big|\,\mathrm{Re}(s)>1\big\rbrace\,\,\mathrm{is\,\,a\,\,compact\,\,subset}.
\end{align*}
Hence, we can rephrase estimate \eqref{woconj} as 
\begin{align*}
Z_{X}(n)=O_{n}\big(Z_{X}(2)\big),
\end{align*}
where implied constant depends on $n$. 

\vspace{0.1cm}
Main Theorem \ref{mainthm1} improves the lower bound \eqref{woconj}, especially the dependence on $n$. However, the upper bound is not optimal in terms of local coordinates, but optimal in terms on $n$.  

\vspace{0.1cm}
The second main result is the following corollary, which is an application of Main Theorem \ref{mainthm1}. 

\vspace{0.1cm}
\begin{maincor}\label{mainthm2}
With notation as above, for any $n\gg1$,  and $\varphi^n\in H^{0}(\mathcal{M}_g,\lambda_n)$ with $X\in \lbrace X_{t}\rbrace$,, we have the following inequality
\begin{align}\label{mainthm2-est}
&\|\varphi^{n}(X)\|_{\mathrm{Qu}}^2 =O_{X_0}\bigg(\frac{(4n^2-4n-3)}{C_{g,1}^{\,6n^2-6n+1}}\bigg(\prod_{j=1}^{j_{0}}\frac{\big|\tau_j(X)\big|^{-\tilde{c}(g,\tau_{j})-2n(n-1)-1/6}}{\big(\log|\tau_{j}(X)|\big)^{6n^2-6n+1}}\bigg)\bigg),
\notag\\[0.1cm]
&\mathrm{where,\,\,for}\,\,1\leq j\leq j_{0},\,\,\tilde{c}(g,\tau_{j}):=e^{-\frac{80(g-1)\log|\tau_j(X)|}{\pi}},
\end{align} 
and implied constant depends on the limiting surface $X_{0}$. Furthermore, $\tau_{j}$ is the local coordinate of $X$, as an element of $\mathcal{M}_{g}$, which is given by the following formula $|\tau_{j}|:=
e^{-2\pi^{2}\slash\ell_{j}}$, and for $1\leq j\leq j_0$, $\ell_{j}$ denotes the length of the closed geodesics, which approach zero, as $X$ approaches $X_{0}$. Moreover, the constant $C_{g,n}$ is as defined in equation \eqref{cgn2}.
\end{maincor}

\vspace{0.1cm}
The third main result is the following corollary, which is an application of Main Theorem \ref{mainthm1} and Corollary \ref{mainthm2}. 

\vspace{0.1cm}
\begin{maincor}\label{mainthm3}
Let $\mathcal{A}\subset\mathcal{M}_g$ be a compact subset. With notation as above, for any $n\geq  2$ and $X\in \lbrace X_{t}\rbrace\cap\mathcal{A}$, we have the following inequality
\begin{align}\label{mainthm3-est}
&\frac{\big|\mu_{g,n}(X)\big|^2}{\prod_{j=1}^{j_0}\big|\tau_j(X)\big|^{n(n-1)}}=O_{\mathcal{A},X_{0}}\big(n^{2}\big),
\end{align} 
and the implied constant depends on the compact subset $\mathcal{A}$, and the limiting surface $X_0$, and $\tau_{j}$  is the local coordinate of $X$, as an element of $\mathcal{M}_{g}$, which is given by the following formula $|\tau_{j}|:=e^{-2\pi^{2}\slash\ell_{j}}$, and for $1\leq j\leq j_0$, $\ell_{j}$ denotes the length of the closed geodesics, which approach zero, as $X$ approaches $X_{0}$.
\end{maincor}

\vspace{0.1cm}
\begin{rem}
For $n=2$, Mumford forms $\mu_{g,2}$ appear in the formula for partition function in Bosonic string theory. Witten for genus one, and D'Hoker and Phong for genus two (series of articles \cite{DP1}--\cite{DP4}) have unified Bosonic and Super string theories. However, the unification is still not complete for genus $g\geq4$. 
\end{rem}

\vspace{0.1cm}
\section{Background material}\label{sec-2}
In this section, we gather results from \cite{MT}, \cite{jj}, \cite{wo}, and \cite{giovanni}, and setup the notation and background material to prove our main results. 

\vspace{0.1cm}
\subsection{Hyperbolic Riemann surface}\label{sec-2.1}
Let
\begin{align*}
\mathbb{H}:=\big\lbrace z=x+iy\in\mathbb{C}\big|\,y=\mathrm{Im}(z)>0\big\rbrace
\end{align*}
denote the hyperbolic upper half-plane, and let $\mu_{\mathrm{hyp}}$ denote the hyperbolic metric on $\mathbb{H}$, which is the natural metric on $\mathbb{H}$, compatible with its complex structure, and is of constant curvature $-1$. Locally, for any $z\in\mathbb{H}$, it is given by the following formula
\begin{align*}
\mu_{\mathrm{hyp}}(z):=\frac{i}{2}\cdot\frac{dz\,d\overline{z}}{(\mathrm{Im}z)^{2}}=\frac{i}{2}\cdot\frac{dx\,dy}{y^{2}},
\end{align*} 
Let $d_{\mathrm{hyp}}(z,w)$ denote the hyperbolic distance between $z=x+iy,w=u+iv$ on $\mathbb{H}$, which is induced by the hyperbolic metric $\mu_{\mathrm{hyp}}$, and is given by the following formula
\begin{align}\label{dhyp}
\cosh^{2}(d_{\mathrm{hyp}}(z,w)\slash 2)=\frac{|z-\overline{w}|^{2}}{4yv}.
\end{align}
Let $X$ denote a compact hyperbolic Riemann surface of genus $g\geq 2$. From Riemann uniformization theorem, $X$ is isometric to $\Gamma_{X}\backslash \mathbb{H}$, where 
$\Gamma_{X}\subset \mathrm{PSL}_{2}(\mathbb{R})$ is a cocompact Fuchsian subgroup acting on the hyperbolic upper half space $\mathbb{H}$, via fractional linear transformations.  Hyperbolic metric descends to $X$, to define a smooth metric on $X$. Locally, for $z,w\in X$, the hyperbolic distance between $z$ and $w$, is given by formula \eqref{dhyp}.

\vspace{0.1cm}
Let
\begin{align*}
\mathrm{vol}_{\mathrm{hyp}}(X):=\int_{X}\mu_{\mathrm{hyp}}(z)=4\pi (g-1)
\end{align*}
denote the hyperbolic volume of $X$.

\vspace{0.1cm}
Let $\mathcal{H}(\Gamma_{X})$ denote the set of inconjugate primitive elements of $\Gamma_{X}$. Let 
\begin{align*}
\ell(\gamma):=\inf_{z\in\mathbb{H}}d_{\mathrm{hyp}}(z,\gamma z),
\end{align*}
denote the length of the closed geodesic determined by $\gamma$ on $X$. Conversely, any closed geodesic on $X$ is determined by a inconjugate primitive element 
$\gamma\in \Gamma_{X}$. Furthermore, for any $n\geq 1$, we have
\begin{align}\label{lgamma-n}
\ell(\gamma^{n})=n\ell(\gamma).
\end{align}
Let $\ell_{X}$ denote the length of the shortest geodesic on $X$.

\vspace{0.1cm}
For any $x\in \mathbb{R}_{>0}$, let 
\begin{align*}
\pi(u):=\big\lbrace  \gamma\in \mathcal{H}(\Gamma_{X})\big|\,e^{\ell_{\gamma}}\leq e^{u}\big\rbrace.
\end{align*}
From Proposition 2.5 on p. 8 in \cite{he}, we have the following estimate 
\begin{align}\label{pgt-1}
\big|\pi(u)\big|=O\big(e^{20D+u}\big),
\end{align}
where $D$ denotes the diameter of $X$, and the implied constant is a universal constant.  

\vspace{0.1cm}
Furthermore, from \cite{chang}, we have the following estimate for $D$, the diameter of $X$
\begin{align}\label{pgt-2}
e^{20D}=O\big(e^{80\pi(g-1)/\ell}\big),
\end{align} 
where $\ell$ denotes the length of the shortest systole on $X$, and the implied constant is a universal constant. 

\vspace{0.1cm}
Hence, combining estimates \eqref{pgt-1} and \eqref{pgt-2}, we arrive at the following primitive version of Prime Geodesic Theorem
\begin{align}\label{pgt}
\big|\pi(u)\big|=O\big(e^{80\pi(g-1)/\ell_{X}+u}\big),
\end{align}
and the implied constant is a universal constant.  

\vspace{0.1cm}
{\textbf{Cotangent bundle.}}
Let $\Omega_{X}^{\otimes n}$ denote the line bundle of holomorphic differentials of degree-$n$, defined over $X$, and let $H^{0}(X,\Omega_{X}^{\otimes n})$ denote the space of holomorphic global sections of $\Omega_{X}^{\otimes n}$. The hyperbolic metric induces a point-wise and $L^2$ metric on $H^{0}(X,\Omega_{X}^{\otimes n})$, denoted by $\|\cdot\|_{\mathrm{hyp}}$ and $\langle\cdot,\cdot\rangle_{\mathrm{hyp}}$, respectively. From Riemann-Roch theorem, we have
\begin{align}
\mathrm{dim}_{\mathbb{C}}\big( H^{0}(X,\Omega_{X}^{\otimes n})\big)=(2n-1)(g-1)+\delta_{1}(n),
\end{align}
where $\delta_{1}(n)$ denotes the Kronecker delta function. 

\vspace{0.1cm}
Furthermore, $\mathcal{S}_{2n}(\Gamma_{X})$, the space of cusp forms of weight-$2n$ with respect to $\Gamma_{X}$, is isometric to $H^{0}(X,\Omega_{X}^{\otimes n})$. At any $z\in X$, locally, any $\alpha\in H^{0}(X,\Omega_{X}^{\otimes n})$ can be represented by $\alpha(z)=f_{\alpha}(z)dz^{\otimes n}$, where $f_{\alpha}\in \mathcal{S}_{2n}(\Gamma_{X})$. Hence, any $\alpha, \,\beta\in H^{0}(X,\Omega_{X}^{\otimes n})$, we have
\begin{align}\label{hyp-met}
\big\langle \alpha,\beta\big\rangle_{\mathrm{hyp}}=\int_{X}y^{2k}f_{\alpha}(z)\overline{f_{\beta}(z)}\mu_{\mathrm{hyp}}(z).
\end{align}

\vspace{0.1cm}
\subsection{Hyperbolic Heat Kernel and Selberg-zeta function}\label{sec-2.2}
At any $z=x+iy\in \mathbb{H}$, let 
\begin{align*}
\Delta_{\mathrm{hyp},0}:=-4y^{2}\frac{\partial^{2}}{\partial z\partial \overline{z}}=-y^{2}\left(\frac{\partial^{2}}{\partial x^2} +\frac{\partial^{2}}{\partial y^2}\right)
\end{align*}
denote the hyperbolic Laplacian acting on smooth functions defined on $\mathbb{H}$. 

\vspace{0.1cm}
Let $K_{\mathbb{H}}:\mathbb{R}_{>0}\times \mathbb{H}\times\mathbb{H}\longrightarrow \mathbb{R}_{>0}$ denote the hyperbolic heat kernel on $\mathbb{H}$, which satisfies the heat equation
\begin{align*}
\bigg(\Delta_{\mathrm{hyp},0}+\frac{\partial }{\partial t}\bigg)K_{\mathbb{H}}(t;z,w)=0,
\end{align*}
where $\Delta_{\mathrm{hyp},0}$ is the hyperbolic Laplacian acting on the $z$-variable.

\vspace{0.1cm} 
An explicit formula for $K_{\mathbb{H}}$ is well-known, and is given by the following formula
\begin{align}\label{kh}
K_{\mathbb{H}}(t;z,w)=\frac{\sqrt{2}e^{-t\slash 4}}{(4\pi t)^{3/2}}\int_{d_{\mathrm{hyp}}(z,w)}^{\infty}\frac{e^{-\rho^{2}\slash 4t}d\rho}{\sqrt{\cosh(\rho)-\cosh(d_{\mathrm{hyp}}(z,w))}}.
\end{align}
The hyperbolic Laplacian $\Delta_{\mathrm{hyp},0}$ descends to define a Laplacian on $X$, which is the natural Laplacian acting on smooth functions defined on $X$. 

\vspace{0.1cm}
The hyperbolic heat kernel $K_{X}:\mathbb{R}_{>0}\times X\times X\longrightarrow\mathbb{R}_{>0}$ is uniquely determined by the partial differential 
\begin{align*}
\bigg(\Delta_{\mathrm{hyp},0}+\frac{\partial }{\partial t}\bigg)K_{X}(t;z,w)=0,
\end{align*}
and the normalization condition
\begin{align*}
\lim_{t\rightarrow 0^{+}}\int_{X}f(w)K_{X}(t;z,w)\mu_{\mathrm{hyp}}(w)=f(z),
\end{align*}
where $f$ is any smooth function on $X$. 

\vspace{0.1cm}
Furthermore, for any $z,w\in X$, we have the following periodization 
\begin{align*}
K_{X}(t;z,w)=\sum_{\gamma\in\Gamma_{X}}K_{\mathbb{H}}(t;z,\gamma w),
\end{align*}
For any $t\in \mathbb{R}_{>0}$ and $z\in X$, set 
\begin{align*}
HK_{X}(t;z):=\sum_{\gamma\in\Gamma_{X}\backslash \lbrace\mathrm{Id}\rbrace}K_{\mathbb{H}}(t;z,\gamma z).
\end{align*}
The hyperbolic heat trace is given by the following equation
\begin{align*}
H\mathrm{Tr}(t):=\int_{X}HK_{X}(t;z)\mu_{\mathrm{hyp}}(z).
\end{align*}
From Selberg trace formula as stated as in \cite{mc}, we have the following formula
\begin{align}\label{sel}
H\mathrm{Tr}K_{X}(t)=\frac{e^{-t\slash 4}}{2\sqrt{4\pi t}}\sum_{n=1}^{\infty}\sum_{\gamma\in\mathcal{H}(\Gamma_{X})}\frac{\ell(\gamma)e^{-\ell^{2}(\gamma^{n})\slash 4t}
}{\sinh(\ell(\gamma^{n})\slash 2)},
\end{align}
where $\mathcal{H}(\Gamma_{X})$ denotes the set of inconjugate primitive elements of $\Gamma_{X}$, and $\ell(\gamma^n)$ denotes the length of the closed geodesic determined by $\gamma^n$ on $X$. 

\vspace{0.1cm}
For any $s\in \mathbb{C}$ with $\mathrm{Re}(s)>1$, the Selberg-zeta function is given by the following formula
\begin{align}\label{selc}
Z_{X}(s):=\prod_{\gamma\in\mathcal{H}(\Gamma_{X})}\prod_{n=1}^{\infty}\big( 1-e^{-(s+k)\ell(\gamma^n)}\big).
\end{align}
From Mckean formula from \cite{mc}, for any $n\geq 1$, we have
\begin{align}\label{mc1}
\frac{Z_{X}^{\prime}(n)}{Z_{X}(n)}=(2n-1)\int_{0}^{\infty}H\mathrm{Tr}K_{X}(t)e^{-n(n-1)t}dt,
\end{align}
where 
\begin{align*}
Z_{X}^{\prime}(n):=\frac{d Z_{X}(s)}{ds}{\bigg\vert_{ s=n}}
\end{align*}
denotes the derivative of $Z_{X}(s)$ at $s=n$.

\vspace{0.1cm}
In connection with the computations associated with the estimates of the Selberg-zeta function, for $t> 2$, we recall the following lower bound from \cite{JK}
\begin{align}\label{jk-sel}
H\mathrm{Tr}K_{X}(t)\geq 1-\mathrm{vol}_{\mathrm{hyp}}(X)K_{\mathbb{H}}(t;0)=1-4\pi(g-1)K_{\mathbb{H}}(t;0).
\end{align} 

\vspace{0.1cm}
\subsection{Hyperbolic Laplacian of degree-$n$}\label{sec-2.3}
Let $X$ be a compact hyperbolic Riemann surface of genus $g\geq 2$, which is isometric to $\Gamma_{X}\backslash\mathbb{H}$. For any $n\geq 1$, at any $z=x+iy\in \mathbb{H}$ (identifying $X$ with its universal cover $\mathbb{H}$), let 
\begin{align*}
\Delta_{\mathrm{hyp},n}:=-y^{2}\bigg(\frac{\partial^{2}}{\partial x^2}+\frac{\partial^{2}}{\partial y^2}\bigg)-2iny\bigg(\frac{\partial }{\partial x}+i\frac{\partial }{\partial y}\bigg)
\end{align*}
denote the hyperbolic Laplacian of degree-$n$, acting on holomorphic differentials of degree-$n$.

\vspace{0.1cm}
As $X$ is compact, $\Delta_{\mathrm{hyp},n}$ admits a discrete spectrum, and let $\lbrace \lambda_{n,k}\rbrace_{k\geq 1}$ denote the set of eigenvalues of 
$\Delta_{\mathrm{hyp},n}$. For any $n\geq 1$, and $s\in \mathbb{C}$ with $\mathrm{Re}(s)>1$, the spectral zeta function associated to the hyperbolic Laplacian $\Delta_{\mathrm{hyp},n}$ is given by the following formula
\begin{align}
\zeta_{X,n}(s):=\sum_{k=1}^{\infty}\lambda_{n,k}^{-s}.
\end{align}
The spectral zeta function admits an analytic continuation to 
\begin{align*}
s\in\mathbb{C}\backslash \big\lbrace  1,1\slash 2-n\big|\,n\in\mathbb{Z}_{\geq 1}\big\rbrace,
\end{align*}
with a pole of order $1$ at $s=1$, and poles of order $2$ at $\lbrace 1\slash 2-n\rbrace_{n\geq 1}$.

\vspace{0.1cm}
Let $\mathrm{det}\Delta_{\mathrm{hyp},n}^{\ast}(X)$ denote the determinant of the regularized determinant, which is given by the following formula
\begin{align*} 
\mathrm{det}\Delta_{\mathrm{hyp},n}^{\ast}(X):=e^{-\zeta_{X,n}^{\prime}(0)},\quad\mathrm{where}\,\,\zeta_{X,n}^{\prime}(0):=\frac{d\zeta_{X,n}(s)}{ds}|_{s=0}.
\end{align*}
From \cite{DP} and Sarnak \cite{sar}, we have
\begin{align}\label{cgn1}
\mathrm{det}\Delta_{\mathrm{hyp},n}^{\ast}(X)=\begin{cases}
C_{g,1}(X)Z_{X}^{\prime}(1)2^{2(g-1)/3+2},\quad\mathrm{for}\,\,n=1;\\
C_{g,n}(X)Z_{X}(n)2^{2(n+1/3)(g-1)},\quad\mathrm{for}\,\,n\geq 2.
\end{cases}
\end{align}
The constant $C_{g,n}$ is explicitly computed in \cite{DP} and \cite{sar}. However, we refer the reader to \cite{giovanni}, as the constant $C_{g,n}$ and the associated computations are described in good detail. Since $X$ is compact, using the isometry $H^{0}(X,\omega^{\otimes n}_{X})\simeq S_{2n}(\Gamma)$, combining Proposition 2.7.2 on p. 63 and Theorem 2.8.4 on p. 77 from \cite{giovanni}, we infer that
\begin{align}\label{cgn2}
&C_{g,n}(X)=e^{-c_{n}\mathrm{vol}_{\mathrm{hyp}}(X)},\notag\\[0.1cm]
&\mathrm{where}\,\,c_{n}:=\frac{\log(G(2n-1))}{2\pi}-\frac{2n-3}{4\pi}\log(\Gamma(2n-1))+\notag\\[0.1cm]&\frac{(2n-1)^{2}}{8\pi}-\frac{(2n-1)\log2\pi}{8\pi}-\frac{\zeta^{\prime}(-1)}{\pi},
\end{align}
and $G(Z)$ denotes the Barnes $G$-function, and $\zeta^{\prime}(-1)$ denotes the derivative of the Riemann-zeta function $\zeta(s)$ at the point $s=-1$.

\vspace{0.1cm}
\subsection{Moduli space of compact Riemann surfaces of genus $g$}\label{sec-2.4}
Let $\mathcal{M}_{g}$ denote the moduli space of isomorphic  classes of compact hyperbolic Riemann surfaces of genus $g\geq 2$. The moduli space $\mathcal{M}_{g}$ is a complex 
orbifold of dimension $3g-3$. Let $\overline{\mathcal{M}}_{g}$ denote the Deligne-Mumford compactification of $\mathcal{M}_{g}$, and the boundary 
\begin{align*}
\partial \overline{\mathcal{M}}_{g}:= \overline{\mathcal{M}}_{g}\backslash\mathcal{M}_{g}=\sum_{j=0}^{[g/2]}\Delta_{j},
\end{align*}
where $\Delta_{j}$ denotes the isomorphism classes of noncompact Riemann surfaces with $j$-nodes. The moduli space $\mathcal{M}_{g}$ admits the structure of a moduli stack, and the boundary $\partial \overline{\mathcal{M}}_{g}$ is a codimenson-$1$ submanifold, and can be realized as a divisor on $\overline{\mathcal{M}}_{g}$. 

\vspace{0.1cm}
Let $\mu_{\mathrm{WP}}$ denote the natural metric on $\mathcal{M}_{g}$, which is not a complete metric.  

\vspace{0.1cm}
We work with a family $\lbrace X_{t}\rbrace$ of compact hyperbolic Riemann surfaces which approach a fixed noncompact hyperbolic Riemann surface 
$X_{0}\in \partial\mathcal{M}_{g}$, as $t$ approaches zero. The limiting Riemann surface $X_{0}$ admits $j_{0}$ punctures, and $k$ connected components.  The family is parametrised by $\mathcal{D}^{n}$, $n$ copies of a slit punctured disc $\mathcal{D}\subset \mathbb{C}$, which is obtained by taking the punctured disc $\mathcal{D}$, and removing a ray connecting the removed origin to a boundary point. 

\vspace{0.1cm}
For $1\leq j \leq j_{0}$, let $\ell_{j}$ denote the length of the closed geodesic on $X_{t}$, which approaches zero, as $t$ approaches zero and $X_{t}$ approaches $X_{0}$. We refer the 
reader to section 2 in \cite{jj}, for further details regarding the family $\lbrace X_{t}\rbrace$.

\vspace{0.1cm}{\textbf{Teichm\"uller modular forms.} }
As described in section \ref{sec-1.2}, let $\pi:\mathcal{C}_{g}\longrightarrow \mathcal{M}_{g}$ be the universal surface. For any $n\geq $1, let $\Lambda_{n}:=\pi_{\ast}(T_{v}\mathcal{C}_{g})^{n}$, where $T_{v}\mathcal{C}_{g}$ denotes the vertical holomorphic tangent bundle of the fibration $\pi$, and the fiber of $\Lambda_{n}$ over any $X\in\mathcal{M}_{g}$ is equal to $H^{0}(X,\Omega_{X}^{\otimes n})$, the space of holomorphic differentials of degree-$n$, defined over the Riemann surface $n$. Let $H^{0}(\mathcal{M}_{g},\lambda_{n})$ denote the space of global holomorphic sections of $\lambda_{n}:=\mathrm{det}(\Lambda_{n})$, the determinant bundle of the vector bundle $\Lambda_{n}$, which are known as Teichm\"uller modular forms. From Riemann-Roch theorem, it follows that the dimension of $H^{0}(\mathcal{M}_{g},\lambda_{n})$, as a complex vector space, is $g_{n}:=(2n-1)(g-1)+\delta_{n}(1)$, where $\delta_{1}(n)$ is the Kronecker delta function. 

\vspace{0.1cm}
Given any $\varphi\in H^{0}(\mathcal{M}_{g},\lambda_{n})$ with trivialization $\varphi:=\varphi_1\wedge\cdots\wedge\varphi_{g_{n}}$, the Quillen metric at any $X\in 
\mathcal{M}_{g}$, is given by the following formula
\begin{align*}
\| \varphi(X)\|_{\mathrm{Qu}}:=\frac{\big|\mathrm{det}N_{\varphi}(X)\big|}{\big|\mathrm{det}^{\ast}\Delta_{\mathrm{hyp},n} (X)\big|},
\end{align*}
where 
\begin{align*}
N_{\varphi}(X):=\left( \big\langle \varphi_i,\varphi_j \rangle_{\mathrm{hyp}}\right)_{1\leq i,j\leq n},
\end{align*}
and the inner-product $\big\langle \cdot,\cdot\rangle_{\mathrm{hyp}}$ is as defined in equation \eqref{hyp-met}.

\vspace{0.1cm}
Let  $H^{0}_{L^{2}, \mathrm{Qu}}(\mathcal{M}_{g},\lambda_{n})\subset H^{0}(\mathcal{M}_{g},\lambda_{n})$ denote the vector subspace space of $L^2$ sections with respect to the $L^2$ Quillen metric. From the behaviour of the Weil-Petersson metric $\mu_{\mathrm{WP}}$ along the Deligne-Mumford boundary, which is as described in \cite{masur}, we infer that
\begin{align*}
H^{0}_{L^2, \mathrm{Qu}}\big(\mathcal{M}_{g},\lambda_{n}\big)\subset H^{0}\big(\overline{\mathcal{M}}_{g},\overline{\lambda}_{n}\otimes O_{\overline{\mathcal{M}}_g}(-\partial\mathcal{M}_g)\big),
\end{align*} 
where $\overline{\lambda}_{n}$ denotes the extension of $\lambda_n$ to $\overline{\mathcal{M}}_g$. 

\vspace{0.1cm}
Furthermore, any $\tilde{\varphi}^n\in H^{0}_{L^{2}, \mathrm{Qu}}(\mathcal{M}_{g},\lambda_{n})$, which is $L^2$ normalized, satisfies the following condition
\begin{align}\label{qu-ip}
\int_{\mathcal{M}_{g}}\|\tilde{\varphi}^n(X)\|_{\mathrm{Qu}}^2\,\mu_{\mathrm{WP}}^{\mathrm{vol}}(X)=1,
\end{align}
where $\mu_{\mathrm{WP}}^{\mathrm{vol}}$ denotes the volume form associated to the Weil-Petersson metric. 

\vspace{0.1cm}
In a seminal article \cite{zograf}, for any $n\geq 1$, $\varphi\in H^{0}(\mathcal{M}_{g},\lambda_{n}) $ and $X\in\mathcal{M}_{g}$, the authors showed that 
\begin{align}\label{c1cal}
c_{1}\big(\lambda_n,\|\cdot\|_{\mathrm{Qu}}\big)(X)=-\frac{i}{2\pi}\partial \overline{\partial}\log \|\varphi\|_{\mathrm{Qu}}=\frac{(6n^2-6n+1)}{3\pi}\mu_{\mathrm{WP}}(X).
\end{align}

\vspace{0.1cm}
{\textbf{Mumford forms.}}
For any $n\geq 1$, from Mumford's isomorphism, we have the following isometry of Hilbert spaces, which are equipped with the Quillen metric 
\begin{align}\label{mumiso}
\mu_{\mathrm{mu}}:  H^{0}(\mathcal{M}_{g},\lambda_{1})^{\otimes (6n^2-6n+1)} \longrightarrow H^{0}(\mathcal{M}_{g},\lambda_{n})
\end{align}

\vspace{0.1cm}
For any $n\geq 1$, the line bundle $\lambda_{n}\otimes\lambda_{1}^{\otimes -(6n^2-6n+1)}$ is trivial, and up to a constant, admits one non-trivial, nowhere vanishing global holomorphic section, which is called the Mumford form, which is  denoted by $\mu_{g,n}$. For $n=2$, the Mumford form $\mu_{g,2}$ appears in the partition function of Bosonic string theory.  

\vspace{0.1cm}
We now recall a special class of global sections of $ H^{0}(\mathcal{M}_{g},\lambda_{n})$ and $ H^{0}(\mathcal{M}_{g},\lambda_{1})$ from \cite{MT}, which will help us in ascertaining the behavior of $\mu_{g,n}$. For any $n\geq 1$, we consider global sections $\varphi^n$, which trivialize at any $X\in\mathcal{M}_{g}$ as $\varphi^n(X):=\varphi_1\wedge\cdots\wedge\varphi_{g_n}$, where $\lbrace \varphi_1,\ldots,\varphi_{g_n}\rbrace$ are a basis of normalized differentials of degree-$n$, which are obtained through Eichler cohomology groups. We refer the reader to section 4 of \cite{MT} for details regarding the construction of the sections $\varphi^n$. 

\vspace{0.1cm}
Furthermore, from Mumford's isometry \eqref{mumiso}, we have 
\begin{align*}
\mu_{\mathrm{mu},n}((\varphi^1)^{\otimes (6n^2-6n+1)})=\varphi^{n},
\end{align*}
from which we derive that, for any $X\in \mathcal{M}_{g}$, we have 
\begin{align}\label{ugn}
\big|\mu_{g,n}(X)\big|^{2}=\bigg|\frac{\mathrm{det}(N_{\varphi^n})(X)}{C_{g,n}Z_{X}(n)}\bigg|\cdot\bigg|\frac{C_{g,1}Z_{X}^{\prime}(1)}{\mathrm{det}(N_{\varphi^1})(X)}\bigg|^{6n^2-6n+1}
\cdot\bigg|\frac{2^{(2/3(g-1)+2)(6n^2-6n+1)}}{2^{2(+1/3)(g-1)}} \bigg|.
\end{align}
From Theorem 1.6 from \cite{fay}, it is known that $\mu_{g,n}$ has a pole of order $n(n-1)/2$ along the boundary $\partial\mathcal{M}_{g}$. 

{\textbf{Asymptotics of Selberg-zeta function.}}
We now recall estimates of Selberg-zeta function from \cite{wo} and \cite{jj}. 

\vspace{0.1cm}
With notation as above, from \cite{wo}, for any $X\in\lbrace X_t\rbrace$, we have the following estimate
\begin{align}\label{zx2}
Z_{X}(2)=O_{X_0}\bigg(\prod_{j=1}^{j_0}\frac{e^{-\pi^2/3\ell_j}}{\ell_{j}^3}\bigg),
\end{align}
where the implied constant depends only on the limiting surface $X_0$.

\vspace{0.1cm}
With notation as above, from \cite{jj}, for any $X\in\lbrace X_t\rbrace$, we have the following estimate
\begin{align}\label{zx1}
\bigg| \frac{Z_{X}^{\prime}(1)}{\mathrm{det}(N_{\varphi^1})(X)} \bigg|=O_{X_{0}}\bigg(\prod_{j=1}^{j_0}\frac{e^{-\pi^2/3\ell_j}}{\ell_{j}}\bigg),
\end{align}
where the implied constant depends only on the limiting surface $X_{0}$. 

\vspace{0.1cm}
With notation as above, for each $1\leq j \leq j_{0}$, let $\tau_{j}$ denote the local coordinate at $X_0$, corresponding to $\ell_j$,  which satisfies the following condition 
\begin{align*}
|\tau_j(X)|:=e^{-2\pi^2/\ell_j}.
\end{align*}
Hence, we can restate estimates \eqref{zx2} and \eqref{zx1} in local coordinates as
\begin{align}\label{zx3}
&Z_{X}(2)=O_{X_0}\bigg(\prod_{j=1}^{j_0}\big(-\log|\tau_{j}(X)|\big)^{3} \big|\tau_j(X)\big|^{1/6}\bigg);\notag\\
&\bigg| \frac{Z_{X}^{\prime}(1)}{\mathrm{det}(N_1)(X)} \bigg|=O_{X_{0}}\bigg(\prod_{j=1}^{j_0}\big(-\log|\tau_{j}(X)|\big)\big|\tau_j(X)\big|^{1/6}\bigg).
\end{align}

\vspace{0.1cm}
Recall that $\mu_{g,n}$ admits a pole of order $n(n-1)/2$ at the boundary $\partial \mathcal{M}_{g}$. Hence, from Theroem 1.6 from \cite{fay}, for any $X\in\lbrace X_t\rbrace$, we have the following estimate
\begin{align}\label{mu-pole}
\big|\mu_{g,n}(X)\big|=O_{X_{0}}\bigg(\prod_{j=1}^{j_0}\big|\tau_j(X) \big|^{-n(n-1)/2}\bigg).
\end{align}
\section{Proofs of Main Results}
In this section, we prove Main Theorem \ref{mainthm1}, and Corollaries \ref{mainthm2}, and \ref{mainthm3}, and infer certain results a corollaries. 
\subsection{Proof of Main Theorem \ref{mainthm1}}\label{sec-3.1}
In this section, we prove Main Theorem \ref{mainthm1}.

\vspace{0.1cm}
\begin{thm}\label{thm1}
With notation as above, for any $n\gg1$, and $X\in \lbrace X_{t}\rbrace$,we have the following estimate
\begin{align}\label{thm1-eqn}
&Z_{X}(2)\ll Z_{X}(n)\ll (4n^2-4n-3)Z_{X}(2)\bigg(\prod_{j=1}^{j_{0}}e^{c(g,\ell_{j})/\ell_{j}^{2}}\bigg),\notag\\&
\mathrm{where,\,\,for}\,\,1\leq j\leq j_{0},\,\,c(g,\ell_{j}):=e^{\frac{160\pi(g-1)}{\ell_j}}.
\end{align} 
\begin{proof}
From equation \eqref{mc1}, for any $s\in \mathbb{C}$ with $\mathrm{Re}(s)>1$, we find
\begin{align}\label{thm1-eqn-1}
\frac{Z_{X}^{\prime}(s)}{Z_{X}(s)}=(2s-1)\int_{0}^{\infty}H\mathrm{Tr}K_{X}(t;0)e^{-s(s-1)t}dt.
\end{align}
Using inequality \eqref{jk-sel}, we arrive at the following inequality
\begin{align}\label{thm1-eqn-2}
\frac{Z_{X}^{\prime}(s)}{Z_{X}(s)}\geq (2s-1)\int_{t_{0}}^{\infty}\big(1-\mathrm{vol}_{\mathrm{hyp}}(X)K_{\mathbb{H}}(t;0)\big)e^{-s(s-1)t}dt.
\end{align}
From formula \eqref{kh}, we choose $t_{g}>2 $ large enough such that, for $t\geq t_{g}$, we have
\begin{align*}
K_{\mathbb{H}}(t;0)\leq \frac{e^{-t\slash 4}}{\mathrm{vol}_{\mathrm{hyp}}(X)}= \frac{e^{-t\slash 4}}{4\pi(g-1)}.
\end{align*}
As the right hand-side of the above inequality is a decreasing function of $g$, the choice of $t_{g}$ is a universal constant satisfying $t_{g}>2$, which we denote by $t_{0}$. So, from inequality \eqref{thm1-eqn-2}, we infer that
\begin{align*}
\frac{Z_{X}^{\prime}(s)}{Z_{X}(s)}\geq (2s-1)\int_{t_0}^{\infty}\big(1-e^{-t/4}\big)e^{-s(s-1)t}dt=(2s-1)\bigg(\frac{e^{-s(s-1)t_0}}{s(s-1)}-\frac{e^{-(s(s-1)+1/4)t_0}}{s(s-1)+1/4}\bigg).
\end{align*}
Integrating, we derive
\begin{align}\label{thm1-eqn-3}
\int_{2}^{n}\frac{Z_{X}^{\prime}(s)ds}{Z_{X}(s)}=\log\big(Z_{X}(n)/Z_{X}(2)\big)\geq \int_{2}^{\infty}(2s-1)\bigg(\frac{e^{-s(s-1)t_0}}{s(s-1)}-\frac{e^{-(s(s-1)+1/4)t_0}}{s(s-1)+1/4}\bigg)ds=\notag\\[0.1cm]\int_{2}^{2+1/4}\frac{e^{-xt_{0}}dx}{x}-\int_{n(n-1)}^{n(n-1)+1/4}\frac{e^{-xt_{0}}dx}{x}=\frac{1}{4}\bigg(\frac{9e^{-9t_{0}/4}}{4}-\frac{e^{-n(n-1)t_0}}{n(n-1)}\bigg)
\end{align}
For $n\gg 1$, it is clear that 
\begin{align}\label{thm1-eqn-4}
\frac{1}{4}\bigg(\frac{9e^{-9t_{0}/4}}{4}-\frac{e^{-n(n-1)t_0}}{n(n-1)}\bigg)\geq c,
\end{align}
where $c>0$ is some positive universal constant.

\vspace{0.1cm}
Combining inequalities \eqref{thm1-eqn-3} and \eqref{thm1-eqn-4}, for $n\gg1$, we arrive at the following inequality
\begin{align}\label{thm1-eqn-5}
Z_{X}(n)\gg Z_{X}(2),
\end{align}  
which proves the lower bound asserted in inequality \eqref{thm1-eqn}.

\vspace{0.1cm}
We now prove the upper bound for $Z_{X}(n)$. Combining equations \eqref{thm1-eqn-1} and \eqref{sel}, we find that
\begin{align*}
\frac{Z_{X}^{\prime}(s)}{Z_{X}(s)}=(2s-1)\int_{0}^{\infty}\frac{e^{-t\slash 4}}{2\sqrt{4\pi t}}\sum_{\gamma\in\mathcal{H}(\Gamma_{X})}\sum_{n=1}^{\infty}\frac{\ell(\gamma)e^{-\ell^{2}(\gamma^{n})\slash 4t}}{\sinh(\ell(\gamma^{n})\slash 2)}e^{-s(s-1)t}dt.
\end{align*}
Replacing the inner-most summation by an integral, and using the fact that $\ell(\gamma^n)=n\ell(\gamma)$, we compute
\begin{align}\label{thm1-eqn-6}
&\frac{Z_{X}^{\prime}(s)}{Z_{X}(s)}\leq\notag\\[0.1cm] &(2s-1)\int_{0}^{\infty}\sum_{\gamma\in\mathcal{H}(\Gamma_{X})}\bigg(\frac{\ell(\gamma)e^{-\ell^{2}(\gamma)\slash 4t}}{\sinh(\ell(\gamma)\slash 2)}+\int_{1}^{\infty}\frac{\ell(\gamma)e^{-\alpha^2\ell^{2}(\gamma)\slash 4t}d\alpha}{\sinh(\alpha\ell(\gamma)\slash 2)} 
\bigg)\frac{e^{-(s(s-1)+1/4)t}}{2\sqrt{4\pi t}}dt.
\end{align}
We now estimate the second term in the summation in the above inequality. For all $x\in \mathbb{R}_{\geq 0}$, using the fact that $\sinh(x)\geq x$, we compute
\begin{align}\label{thm1-eqn-7}
\frac{1}{\sqrt{4\pi t}}\int_{1}^{\infty}\frac{\ell(\gamma)e^{-\alpha^2\ell^{2}(\gamma)\slash 4t}d\alpha}{\sinh(\alpha\ell(\gamma)\slash 2)}=
\frac{1}{\sqrt{4\pi t}}\int_{1}^{\infty}2\alpha e^{-\alpha^2\ell^{2}(\gamma)\slash 4t}d\alpha\ll\notag\\[0.1cm]
\sqrt{t}\int_{\ell^{2}(\gamma)/4t}^{\infty}\frac{e^{-\theta}d\theta}{\ell^{2}(\gamma)}=\frac{\sqrt{t}e^{-\ell^{2}(\gamma)/4t}}{\ell^{2}(\gamma)}.
\end{align}
Using counting function \eqref{pgt}, we now derive
\begin{align}\label{thm1-eqn-8}
\sum_{\gamma\in\mathcal{H}(\Gamma_{X})}\frac{\sqrt{t}e^{-\ell^{2}(\gamma)/4t}}{\ell^{2}(\gamma)}\ll
\bigg(\sum_{j=1}^{j_{0}}e^{\frac{160\pi(g-1)}{\ell_j}}\frac{\sqrt{t}e^{-\ell^{2}_{j}/4t}}{\ell^{2}_{j}}+ \sqrt{t}\sum_{j=1}^{j_{0}}e^{\frac{160\pi(g-1)}{\ell_j}}
\int_{\ell_{j}}^{\infty}\frac{e^{-u^{2}/4t+u}du}{u^{2}}\bigg),
\end{align}
where, for $1\leq j\leq j_{0}$, $\ell_{j}$ denotes the length of the closed geodesics on $X$, which approaches zero, as $t$ approaches zero. 

\vspace{0.1cm}
We now compute
\begin{align}\label{thm1-eqn-9}
\sqrt{t}\sum_{j=1}^{j_{0}}\int_{\ell_{j}}^{\infty}\frac{e^{-u^{2}/4t+u}du}{u^{2}}\leq \sum_{j=1}^{j_{0}}\frac{\sqrt{t}e^{t}}{\ell^{2}_{j}}\int_{-\infty}^{\infty}e^{-(u/2\sqrt{t}-\sqrt{t})^2}du\ll
 \sum_{j=1}^{j_{0}}\frac{te^{t}}{\ell^{2}_{j}}
\end{align}
Combining inequalities \eqref{thm1-eqn-7}--\eqref{thm1-eqn-9}, we arrive at the following inequality
\begin{align}\label{thm1-eqn-10}
\frac{1}{\sqrt{4\pi t}}\sum_{\gamma\in\mathcal{H}(\Gamma_{X})}\int_{1}^{\infty}\frac{\ell(\gamma)e^{-\alpha^2\ell^{2}(\gamma)\slash 4t}d\alpha}{\sinh(\alpha\ell(\gamma)\slash 2)} 
\ll  \sum_{j=1}^{j_{0}}\frac{te^{\frac{160\pi(g-1)}{\ell_j}+t}}{\ell^{2}_{j}}.
\end{align}
We now estimate the first term in the summation in inequality \eqref{thm1-eqn-6}. From similar arguments as the ones used to prove inequalities \eqref{thm1-eqn-7}--\eqref{thm1-eqn-9}, we derive
\begin{align}\label{thm1-eqn-11}
\frac{1}{\sqrt{4\pi t}}\sum_{\gamma\in\mathcal{H}(\Gamma_{X})}\frac{\ell(\gamma)e^{-\ell^{2}(\gamma)\slash 4t}}{\sinh(\ell(\gamma)\slash 2)}\ll\frac{1}{\sqrt{4\pi t}}
\sum_{j=1}^{j_{0}}e^{-\ell^{2}_{j}\slash 4t}+\frac{1}{\sqrt{4\pi t}}\sum_{j=1}^{j_{0}}e^{\frac{160\pi(g-1)}{\ell_j}}\int_{\ell_{j}}^{\infty}e^{-u^2/4t+u}du\ll\notag\\
\frac{1}{\sqrt{4\pi t}}\sum_{j=1}^{j_{0}}e^{-\ell^{2}_{j}\slash 4t}+ \sum_{j=1}^{j_{0}}e^{\frac{160\pi(g-1)}{\ell_j}+t}.
\end{align}
From elementary calculus, for any $1\leq j\leq j_{0}$, it follows that
\begin{align}\label{thm1-eqn-12}
\max_{t\in\mathbb{R}_{\geq 0}}\frac{1}{\sqrt{4\pi t}}e^{-\ell^{2}_{j}\slash 4t}=\bigg(\frac{1}{\sqrt{4\pi t}}e^{-\ell^{2}_{j}\slash 4t}\bigg)_{t=\ell_{j}/2}\ll\frac{1}{\sqrt{\ell_j }}.
\end{align}
Combining inequalities \eqref{thm1-eqn-10}--\eqref{thm1-eqn-12}, we compute
\begin{align*}
&\int_{0}^{\infty}\frac{e^{-(s(s-1)+1/4)t}}{2\sqrt{4\pi t}}\sum_{\gamma\in\mathcal{H}(\Gamma_{X})}\bigg(\frac{\ell(\gamma)e^{-\ell^{2}(\gamma)\slash 4t}}{\sinh(\ell(\gamma)\slash 2)}+\int_{1}^{\infty}\frac{\ell(\gamma)e^{-\alpha^2\ell^{2}(\gamma)\slash 4t}d\alpha}{\sinh(\alpha\ell(\gamma)\slash 2)} \bigg)\ll\notag\\[0.1cm]
&\sum_{j=1}^{j_{0}}\frac{1}{\sqrt{\ell_j}}\int_{0}^{\infty}e^{-(s(s-1)+1/4)t}dt+\sum_{j=1}^{j_{0}}e^{\frac{160\pi(g-1)}{\ell_j}}\bigg(\int_{0}^{\infty}e^{-(s(s-1)-1/4)t}dt+\int_{0}^{\infty}\frac{te^{-(s(s-1)-3/4)t}dt}{\ell^{2}_{j}}\bigg).
\end{align*}
Assuming that $s\in\mathbb{R}_{\geq 2}$, and applying the fact that, for $1\leq j\leq j_{0}$, $\ell_{j}$ approaches zero, from the above inequality, we infer that
\begin{align}\label{thm1-eqn-13}
\int_{0}^{\infty}\sum_{\gamma\in\mathcal{H}(\Gamma_{X})}\bigg(\frac{\ell(\gamma)e^{-\ell^{2}(\gamma)\slash 4t}}{\sinh(\ell(\gamma)\slash 2)}+\int_{1}^{\infty}\frac{\ell(\gamma)e^{-\alpha^2\ell^{2}(\gamma)\slash 4t}d\alpha}{\sinh(\alpha\ell(\gamma)\slash 2)} \bigg)\frac{e^{-(s(s-1)+1/4)t}}{2\sqrt{4\pi t}}dt\ll\notag\\\frac{e^{\frac{160\pi(g-1)}{\ell_j}}}{s(s-1)-3/4}\sum_{j=1}^{j_0}\frac{1}{\ell_{j}^2}.
\end{align}
Combining inequalities \eqref{thm1-eqn-6} and \eqref{thm1-eqn-13}, and integrating, we derive
\begin{align*}
\int_{2}^{n}\frac{Z_{X}^{\prime}(s)}{Z_{X}(s)}=\log\big(Z_{X}(n)/Z_{X}(2)\big)\ll \bigg(\sum_{j=1}^{j_0}\frac{e^{\frac{160\pi(g-1)}{\ell_j}}}{\ell_{j}^2}\bigg)\int_{2}^{n}\frac{(2s-1)ds}{s(s-1)-3/4}
= \\[0.1cm]\bigg(\sum_{j=1}^{j_0}\frac{e^{\frac{160\pi(g-1)}{\ell_j}}}{\ell_{j}^2}\bigg)\log\big((4n^{2}-4n-3)/5\big),
\end{align*}
which implies that
\begin{align*}
&Z_{X}(n)\ll (4n^2-4n-3)Z_{X}(2)\bigg(\prod_{j=1}^{j_{0}}e^{c(g,\ell_{j})/\ell_{j}^{2}}\bigg),\\&
\mathrm{where,\,\,for}\,\,1\leq j\leq j_{0},\,\,c(g,\ell_{j}):=e^{\frac{160\pi(g-1)}{\ell_j}},
\end{align*}
and proves the upper-bound asserted in \eqref{thm1-eqn}, and completes the proof of the theorem.
\end{proof}
\end{thm}

\vspace{0.1cm}
\begin{cor}\label{cor2}
With notation as above, for any $n\gg1$, and $X\in \lbrace X_{t}\rbrace$, we have the following inequalities
\begin{align}\label{cor2-eqn}
&\prod_{j=1}^{j_{0}}\big(-\log|\tau_j(X)|\big)^{3}\big|\tau_{j}(X)\big|^{1/6}\ll_{X_{0}} Z_{X}(n)\ll_{X_{0}} (4n^2-4n-3)\bigg(\prod_{j=1}^{j_{0}}\frac{\big(-\log|\tau_j(X)|\big)^{3}}{\big|
\tau_{j}(X)\big|^{\tilde{c}(g,\tau_{j})-1/6}}\bigg),\notag\\&
\mathrm{where,\,\,for}\,\,1\leq j\leq j_{0},\,\,\tilde{c}(g,\tau_{j}):=e^{-\frac{80(g-1)\log|\tau_j(X)|}{\pi}},
\end{align} 
and the implied constants in both the inequalities, depend only on the limiting surface $X_0$.
\begin{proof}
Proof of the corollary follows from combining inequality \eqref{thm1-eqn} with equation \eqref{zx3}.
\end{proof}
\end{cor}
\subsection{Proof of Corollaries \ref{mainthm2} and \ref{mainthm3}}\label{sec-3.2}
 In this section, using estimate \eqref{thm1-eqn}, we prove Corollaries \ref{mainthm2} and \ref{mainthm3}.

\vspace{0.1cm}
\begin{prop}\label{prop3}
With notation as above, for any $n\gg1$,  and $\varphi^n\in H^{0}(\mathcal{M}_g,\lambda_n)$ with $X\in \lbrace X_{t}\rbrace$, as $t$ approaches zero, we have the following inequalities
\begin{align}\label{prop3-eqn}
&\frac{e^{cgn^2}C_{g,n}}{C_{g,1}^{\,6n^2-6n+1}}\bigg(\prod_{j=1}^{j_{0}}\frac{|\tau_{j}(X)|^{-2n(n-1)}}{\big(\log|\tau_{j}(X)|\big)^{6n^2-6n-2}}\bigg)\ll_{X_{0}}  \big|\mathrm{det}(N_{\varphi^n})(X)\big| \ll_{X_{0}}\notag\\& \frac{(4n^2-4n-3)e^{cgn^2}C_{g,n}}{C_{g,1}^{\,6n^2-6n+1}}\bigg(\prod_{j=1}^{j_{0}}\frac{\big|\tau_j(X)\big|^{-\tilde{c}(g,\tau_{j})-2n(n-1)}}{\big(\log|\tau_{j}(X)|\big)^{6n^2-6n-2}}\bigg),
\notag\\
&\mathrm{where,\,\,for}\,\,1\leq j\leq j_{0},\,\,\tilde{c}(g,\tau_{j}):=e^{-\frac{80(g-1)\log|\tau_j(X)|}{\pi}},
\end{align} 
and $c>0$ is a universal constant, and the implied constants in both the inequalities, depend only on the limiting surface $X_0$.
\begin{proof}
From Mumford isomorphism \eqref{mumiso} and combining equations \eqref{ugn} and \eqref{mu-pole}, as $t$ approaches zero, we have 
\begin{align}\label{prop3-eqn1}
\big|\mu_{g,n}(X)\big|^2=\bigg|\frac{\mathrm{det}(N_{\varphi^n})(X)}{C_{g,n}Z_{X}(n)}\bigg|\cdot\bigg|\frac{C_{g,1}Z_{X}^{\prime}(1)}{\mathrm{det}(N_{\varphi^1})(X)}\bigg|^{6n^2-6n+1}\cdot\bigg|\frac{(2^{2(g-1)/3+2)(6n^-6n+1)}}{2^{2(n+1/3)(g-1)}}\bigg|=
\notag\\ O_{X_{0}}\bigg(e^{cgn^2}\prod_{j=1}^{j_0}\big|\tau_j(X) \big|^{-n(n-1)/2}\bigg),
\end{align}
where $c>0$ is a universal constant.

\vspace{0.1cm}
Proof of the proposition follows from combining equations \eqref{prop3-eqn1} and \eqref{zx3} with inequalities described in equation \eqref{cor2-eqn} completes the proof of the theorem.
\end{proof}
\end{prop}

\vspace{0.1cm}
\begin{cor}\label{cor4}
With notation as above, for any $n\gg1$,  and $\varphi^n\in H^{0}(\mathcal{M}_g,\lambda_n)$ with $X\in \lbrace X_{t}\rbrace$, we have the following estimate
\begin{align*}
&\|\varphi^n(X)\|_{\mathrm{Qu}}^2=O_{X_0}\bigg((4n^2-4n-3)e^{cgn^2}\bigg(\prod_{j=1}^{j_{0}}\frac{\big|\tau_j(X)\big|^{-\tilde{c}(g,\tau_{j})-2n(n-1)-1/6}}{\big(\log|\tau_{j}(X)|\big)^{6n^2-6n+1}}\bigg)\bigg),
\notag\\[0.1cm]
&\mathrm{where,\,\,for}\,\,1\leq j\leq j_{0},\,\,\tilde{c}(g,\tau_{j}):=e^{-\frac{80(g-1)\log|\tau_j(X)|}{\pi}},
\end{align*} 
and $c>0$ is a universal constant, and the implied constant depends only on the limiting surface $X_0$.
\begin{proof}
From the definition of the Quillen metric, we have
\begin{align}\label{cor4-eqn1}
\|\varphi^n(X)\|_{\mathrm{Qu}}^2=\bigg|\frac{\mathrm{det}(N_{\varphi^n})(X)}{2^{2(n+1/3)(g-1)}C_{g,n}Z_{X}(n)}\bigg|.
\end{align}
The proof of the corollary follows from combining Theorem \ref{thm1} and Proposition \ref{prop3} with equation \eqref{cor4-eqn1}.
\end{proof}
\end{cor}

\vspace{0.1cm}
\begin{cor}\label{cor5}
Let $\mathcal{A}\subset \mathcal{M}_{g}$ be a compact subset. With notation as above, for any $n\gg1$,  and $\varphi^n\in H^{0}(\mathcal{M}_g,\lambda_n)$ with $X\in \lbrace X_{t}\rbrace$, we have the following estimate
\begin{align*}
\|\varphi^n(X)\|_{\mathrm{Qu}}^2=O_{\mathcal{A},X_0}\big(n^{2}e^{cgn^2} \big),
\end{align*}
where $c>0$ is a universal constant, and the implied constant depends on the compact subset $\mathcal{A}\subset \mathcal{M}_g$, and the limiting surface $X_0$.
\begin{proof}
The proof of the corollary follows from Corollary \eqref{cor4}.
\end{proof}
\end{cor}

\vspace{0.1cm}
\begin{cor}\label{cor6}
Let $\mathcal{A}\subset \mathcal{M}_{g}$ be a compact subset. With notation as above, for $X\in \lbrace X_{t}\rbrace\cap \mathcal{A}$, and $n\gg1$, we have the following estimate
\begin{align*}
&\frac{\big|\mu_{g,n}(X)\big|^2}{\prod_{j=1}^{j_0}\big|\tau_j(X)\big|^{n(n-1)}}=O_{\mathcal{A},X_{0}}\big(n^{2}\big),
\end{align*} 
and the implied constant depends on the compact subset $\mathcal{A}$, and the limiting surface $X_0$.
\begin{proof}
The proof of the corollary follows from  combining equation \eqref{ugn} with estimates \eqref{zx3}, \eqref{thm1-eqn}, \eqref{prop3-eqn}.
\end{proof}
\end{cor}

\vspace{0.1cm}
\begin{thm}\label{thm7}
With notation as above, for $n\gg1$, let $\tilde{\varphi}^n\in H^{0}_{L^2,\mathrm{Qu}}(\mathcal{M}_g),\lambda_{n})$ denote an $L^2$ normalized section, i.e., $\tilde{\varphi}^n$ satisfies equation \eqref{qu-ip}. Then, for any $X\in\mathcal{M}_g$, we have the following estimate  
\begin{align*}
\limsup_{n\rightarrow \infty}\|\tilde{\varphi}^n(X)\|_{\mathrm{Qu}}^{2}=O\bigg(\bigg(\frac{6n^2-6n+1)}{3\pi}\bigg)^{3g-3}\bigg);
\end{align*}
and for $\mathcal{A}\subset\mathcal{M}_g$ compact subset, we have the following estimate
\begin{align}\label{thm7-eqn}
\|\tilde{\varphi}^n(X)\|_{\mathrm{Qu}}^{2}=O_{\mathcal{A}}\bigg(\bigg(\frac{6n^2-6n+1)}{3\pi}\bigg)^{3g-3}\bigg),
\end{align}
where the implied constant depends on the compact subset $\mathcal{A}\subset\mathcal{M}_g$.
\begin{proof}
Let $(M,\omega)$ be a complex manifold of dimension $r$, and let $\ell$ be a holomorphic line defined over $M$. For any $n\geq 1$, let $H^0(M,\ell^{\otimes n})$ denote the vector space of global holomorphic sections of the line bundle $\ell^{\otimes n}$. Let $\|\cdot\|_{\ell^{\otimes n}}$ and $\langle\cdot,\cdot\rangle_{\ell^{\otimes n}}$ denote point-wise and $L^2$ norms on  $H^0(M,\ell^{\otimes n})$, respectively.  We are assuming that $H^0(X,\ell^{\otimes n})$ is finite dimensional. 

\vspace{0.1cm}
In normal coordinates, for any $s\in H^{0}(M,\ell)$, at any $z\in M$, we have
\begin{align*}
&\omega(z):=\frac{i}{2}\sum_{j=1}^{r} dz_j\wedge d\overline{z_j},\\[0.1cm]
&c_1(\ell,\|s\|_{\ell})(z):=-\frac{i}{2\pi}\partial\overline{\partial}\log\|s\|_{\ell}^2=\frac{i}{2}\sum_{j=1}^{r}\lambda_{j}dz_j\wedge d\overline{z_j}.
\end{align*}
Let $s\in H^{0}(M,\ell^{\otimes n})$ be an $L^2$ normalized section. Then, in \cite{berman} proved the following estimate
\begin{align}\label{thm7-eqn1}
\limsup_{k\rightarrow \infty}\|s(z)\|_{\ell^{\otimes n}}^2=O\bigg(\big(\prod_{j=1}^{r}\lambda_j\big) n^r\bigg).
\end{align}
We now apply Berman's estimate to the setting of $L^2$ normalized Teichm\"uller modular forms, i.e., we set $X=\mathcal{M}_g$, $\omega=\mu_{\mathrm{WP}}$, $n=3g-3$, and 
$\ell_1=\lambda_1$. Recall that, from Mumford's isomorphism,  we have $\lambda_n=\lambda_1^{\otimes (6n^2-6n+1)}$, and from Zograf-Taktajan's result \eqref{c1cal}, for any $\tilde{\varphi}^{1}\in H^{0}(\mathcal{M}_g,\lambda_1)$, at $X\in\mathcal{M}_g$, we find
\begin{align}\label{thm7-eqn2}
c_1(\lambda_1,\|\tilde{\varphi}^{1}\|_{\mathrm{Qu}})(X)=\frac{1}{3\pi}\mu_{\mathrm{WP}}(X). 
\end{align}
Hence, from equation \eqref{thm7-eqn2}, and applying estimate \eqref{thm7-eqn1} to our setting, for any $L^2$ normalized section $\tilde{\varphi}^n\in H^{0}_{L^2,\mathrm{Qu}}(\mathcal{M}_g,\lambda_{n})$, we have the following estimate 
\begin{align}\label{thm7-eqn3}
\limsup_{n\rightarrow \infty}\|\tilde{\varphi}^n(X)\|_{\mathrm{Qu}}^2=O\bigg(\bigg(\frac{6n^2-6n+1)}{3\pi}\bigg)^{3g-3}\bigg).
\end{align}
Estimate \eqref{thm7-eqn} follows directly follows from estimate \eqref{thm7-eqn3}.
\end{proof}
\end{thm}

\vspace{0.2cm}
\subsection*{Acknowledgements} Both the authors thank Prof. Nagaraj, Prof. Biswas, and Dr. Morye for valuable discussions.   

\end{document}